\documentclass[12pt]{article}
\usepackage[english]{babel}
\usepackage[cp1251]{inputenc}
\usepackage{amssymb}
\usepackage{amsthm}
\usepackage{inputenc}
\begin{document}
\newtheorem{defn}{Definition}[section]
\newtheorem{lm}{Lemma}[section]
\newtheorem{thm}{Theorem}[section]
\newtheorem{pr}{Proposition}[section]
\newtheorem{exam}{Example}[section]
\newtheorem{cor}{Corollary}[section]
\author{Sh.A. Ayupov  $^{1,*},$  R.Z. Abdullaev $^2$, K.K. Kudaybergenov $^3$}

\title{\bf On a certain class of operator algebras and their derivations}

\maketitle
\begin{abstract}

 Given a von Neumann algebra $M$ with a faithful normal finite trace, we introduce the so called finite tracial
 algebra $M_f$ as the intersection of $L_p$-spaces $L_p(M, \mu)$ over all $p\geq1$ and over all faithful normal finite
 traces $\mu$ on $M.$ Basic algebraic and topological properties of finite tracial algebras are studied. We prove that all derivations on these algebras are inner.
\end{abstract}

\medskip

$^1$ \emph{Institute of Mathematics and Information Technologies, Uzbekistan Academy
 of Science, Dormon Yoli str. 29, 100125, Tashkent, Uzbekistan}

 \emph{and}

 \emph{Abdus Salam International Centre for Theoretical Physics, Trieste, Italy,}

e-mail: \emph{sh\_ayupov@mail.ru}

 $^{2}$ \emph{Institute of Mathematics and Information Technologies, Uzbekistan
Academy of Science, Dormon Yoli str. 29, 100125, Tashkent, Uzbekistan,}

e-mail: \emph{arustambay@yandex.ru}

 $^{3}$ \emph{Karakalpak state university, Ch. Abdirov str. 1, 142012, Nukus, Uzbekistan,}

e-mail: \emph{karim2006@mail.ru}



\medskip
* Corresponding author

\newpage
\section{Introduction}

In the present paper we introduce a new class of algebras, the  so called \emph{finite tracial
algebras}, which are defined as the intersection of non commutative $L_p$-spaces $L_p(M,\mu)$
\cite{Yed} over all $p\in [1,\infty)$  and over all faithful normal finite (f.n.f.) traces
$\mu$ on a von Neumann  algebra $M.$ Equivalently,  a finite tracial algebra $M_f$ is the
intersection of all non commutative Arens algebras  $L^\omega(M,\mu)=\bigcap\limits_{p\geq 1}L_p(M,\mu),$
over all f.n.f. traces $\mu.$ It is known that Arens algebras are metrizable locally convex *-algebras
with respect to the topology generated by the system of $L_p$-norms for a fixed trace. Algebraic and
 topological properties of Arens algebras have been investigated
 in the papers \cite{Abd}- \cite{Alb},  \cite{Are}, \cite{Ino}.

In the present paper we study basic properties of finite tracial algebras with the topology generated by
all $L_p$-norms $\{\|\cdot\|^{\mu}_{p}\}$, where $p\in [1,\infty)$  and $\mu$  runs over all f.n.f. traces
on the given von Neumann algebra $M.$ We prove that a finite tracial algebra $M_f$ is metrizable or reflexive
if and only if the center  of the von Neumann algebra $M$ is finite dimensional; in this case $M_f$
coincides with an appropriate Arens algebra. We also give a necessary and sufficient condition for $M_f$
to coincide (as a set) with $M.$ But even in this case one has a new topology on the  von Neumann algebra $M.$
We obtain also a description of the dual space for the algebra $M_f.$

Finally we prove that every derivation on a solid subalgebra of the Arens algebra $L^\omega(M,\tau)$ is inner.
In particular we obtain that the algebra $M_f$ admits only inner derivations.

Throughout the paper we consider a von Neumann algebra $M$ with a f.n.f trace.
Therefore $M$ is a finite
 von Neumann algebra and thus all closed densely defined operators affiliated with $M$
  are measurable with respect to $M$,
 i.~e. the set of all such operators coincides with the algebra $S(M)$ of all
 measurable operators and hence also with the algebra $LS(M)$
 of all locally measurable operators affiliated with $M$ ; moreover the center of
  $S(M)=LS(M)$ coincides with the set of operators
 affiliated  with the center
 of $M.$

 \section{Preliminaries}

Let $M$ be  a von Neumann algebra with the positive
cone $M^+$ and let $\textbf{1}$ denote the identity operator in $M.$

A positive linear functional $\mu$ is called a \emph{finite trace} if  $\mu(u^\ast xu)=\mu(x)$
for all $x\in M$ and each unitary operator $u\in M.$

A finite trace $\mu$ is said to be \emph{faithful} if  for $x\in M^{+},$ $\mu(x)=0$ implies that $x=0.$

A finite trace $\mu$ is \emph{normal} if given any monotone net $\{x_\alpha\}$ increasing to $x\in M$, one
has $\mu(x)=\sup\mu(x_\alpha).$

Let $\tau$ be a fixed faithful normal finite  (f.n.f.) trace on a von Neumann algebra $M.$
The Radon~--- Nikodym theorem  \cite[Theorem 14]{Seg} implies that given any f.n.f. trace $\mu$ on $M$
there exists a positive operator $h\in L^{1}(M,\tau)$ affiliated  with the center of $M$
 such that $\mu(x)=\tau(hx)$ for all $x\in M.$
This operator $h$ is called the  Radon~--- Nikodym derivative of the trace $\mu$ with respect to
the trace $\tau$ and it is denoted as $\frac{\textstyle d\mu}{\textstyle d\tau}.$

We recall \cite{Seg}, \cite{Yed} that given a f.n.f. trace $\tau$ on a von Neumann algebra $M$ the space
$L_p(M,\tau)$ , $p\in [1,\infty)$, is defined as
$$L_p(M,\tau)=\{x\in S(M):\,\,|x|^p\in L_1(M,\tau)\}.$$
The space $L_p(M,\tau)$ equipped with the norm $\|x\|_p=(\tau(|x|^p))^{1/p}$ is a Banach space and
its dual space coincides with $L_q(M,\tau)$ where
$\frac{\textstyle 1}{\textstyle p}+\frac{\textstyle 1}{\textstyle q}=1,$ and the
duality is given by
$$\langle x,a\rangle=f_a(x)=\tau(ax),$$
for all $f_a\in L_{p}(M,\tau)^{\ast},\,\,\,\,a\in L_{q}(M,\tau)$  (see \cite[Theorem 4.4]{Yed}).

Following \cite{Ino} consider the intersection
$$L^\omega(M,\tau)=\bigcap\limits_{p\in [1,\infty)}L_p(M,\tau).$$
It is known (see also \cite{Abd}, \cite{Alb}, \cite{Are}),  that $L^\omega(M,\tau)$
is a complete locally convex $\ast$-algebra with  respect to the topology $t^{\tau}$ generated by the system of
norms  $\{\|\cdot\|_p\}_{p\in[1,\infty)}.$

Each operator $a \in\bigcup\limits_{q\in(1,\infty)}L_q(M,\tau)$ defines a continuous linear
functional $f_a$ on $(L^\omega(M,\tau),t^\tau)$ by the formula $f_a(x)=\tau (ax),$ and conversely
given an arbitrary continuous linear functional $f$ on the algebra  $(L^\omega(M,\tau),t^\tau)$
there exists an element  $a\in\bigcup\limits_{q\in(1,\infty)}L_q(M,\tau)$
such that $f(x)=\tau (ax).$

 \section{Finite Tracial Algebras}

 Let $M$ be  a finite von Neumann algebra. Denote by $\mathcal{F}$ the set of all f.n.f. traces on $M$
 and from now on suppose that $\mathcal{F}\neq\varnothing.$

 Consider the space
$$M_f=\bigcap\limits_{\mu\in \mathcal{F}}\bigcap\limits_{p\in[1,\infty)}
L_{p}(M,\mu)=\bigcap_{\mu\in\mathcal{F}}L^{\omega}(M, \mu).$$
On the space  $M_f$ one can consider the topology $t,$ generated by the system  of norms
$\{\|\cdot\|_{p}^{\mu}: \mu\in \mathcal{F}, p\in[1,\infty)\}.$

Since each Arens algebra $L^{\omega}(M,\mu),\,\,\mu\in\mathcal{F},$ is a complete locally convex
topological $\ast$-algebra in $S(M)$ from the above definition one easily obtains the following

\begin{thm}
$(M_f,t)$ is a complete locally convex topological $\ast$-algebra.
\end{thm}

\textbf{Definition.}
The topological $\ast$-algebra $M_f$ is called the  \emph{finite tracial algebra} with respect to the
von Neumann algebra $M.$

\textbf{Remark.} Finite tracial algebras present examples of so called $GW^{\ast}$-algebras in the sense of
\cite{Kun}.

Recall (see \cite{Kun}) that a topological $\ast$-algebra $(A, t_A)$ is called
$GW^{\ast}$-algebra, if $A$ has a
$W^{\ast}$-subalgebra  $B$ with $(\textbf{1}+x^{\ast}x)^{-1}\in B$ for all
$x\in A$ and the unit ball of $B$ if $t_A$-bounded.

The finite tracial algebra $M_f$ is a $GW^{\ast}$-algebra. Since $M\subset M_f$ it
is sufficient to show that the unit ball in $M$ is $t$ - bounded in $M_f$.

Let $x\in M$, $\|x\|_\infty\leq 1$. For $\mu\in \mathcal{F}$, and $1\leq p< \infty$ we have
$$\|x\|_{p}^{\mu}=\|x\textbf{1}\|_p^\mu\leq \|x\|_\infty\|\textbf{1}\|_p^\mu\leq \mu(\textbf{1})^{\frac{1}{p}},$$
i. e. $\|x\|_p^\mu\leq \mu(\textbf{1})^{\frac{1}{p}}$ for all $x\in M$, $\|x\|_\infty\leq 1$.
 This means that the unit ball of $M$
is $t$ - bounded in $M_f$. Therefore $M_f$ is a $GW^*$ - algebra.

The algebra $M_f$ contains $M$ but it is a rather small algebra, since it is contained in all $L_p(M,\mu)$
for all $p\geq 1$ and f.n.f. traces $\mu$ on $M.$ The following result gives necessary and sufficient
conditions for $M_f$ to coincide with $M.$

\begin{thm}
For a finite von Neumann algebra $M$ the following conditions are equivalent

i) $M_f=M;$

ii) $M$ is a finite sum of homogeneous type $I_n,\, n\in\mathbb{N}$ von Neumann algebras.
\end{thm}

The proof  of this theorem consists of several auxiliary proposition which are interesting
in their own right. Let us start with the commutative case.

\begin{pr}
Let $M$ be a von Neumann algebra with a faithful normal trace and let $Z$ be its center. Then the center
of the algebra $M_f$ coincides with $Z$, i. e. $Z(M_f)=Z$. In particular if $M$ is commutative then $M_f=M$.
\end{pr}

\emph{Proof.} Let $M$ be a von Neumann algebra with a faithful normal finite trace $\tau$, and $\tau(\textbf{1})=1$.

Consider $x\in Z(M_f)$, $x\geq 0$, and let $x=\int\limits_0^\infty \lambda d e_\lambda$
be the spectral resolution of $x$.
Since $x\in Z(M_f)$ and $M\subset M_f$, we have that $e_\lambda \in Z$
for all $\lambda\in \mathbb{R}$. Passing if necessary
to the element $\varepsilon \textbf{1}+x$ we may suppose without loss of generality that $e_1=0$.

For $n\in \mathbb{N}$ set
$$p_n=e_{(n+1)^2}-e_{n^2}$$
and
$$y=\sum_{n\in \mathbb{N}}n^2 p_n.$$
Since $xp_n\geq n^2p_n$ for all $n\in \mathbb{N}$, we have that $0\leq y\leq x$ and hence $y\in M_f$.

Let
$$F=\{n\in \mathbb{N}:t_n=\tau(p_n)\neq 0\}$$
and
$$h=\sum_{n\in F}\frac{1}{n^2 t_n} p_n \in Z(S(M)).$$
Since
$$\bigvee\limits_{n=1}^{m}p_n=\bigvee\limits_{n=1}^{m}(e_{(n+1)^2}-e_{n^2})=
\sum\limits_{n=1}^{m}(e_{(n+1)^2}-e_{n^2})=e_{(m+1)^2}-e_1=e_{(m+1)^2}\uparrow \textbf{1},$$
one has that
$$\bigvee\limits_{n=1}^{\infty} p_n=\textbf{1}.$$
Therefore there exists $h^{-1}\in S(M)$. Further we have
$$\tau(h)=\sum_{n\in F}\frac{1}{n^2t_n}\tau(p_n)=\sum_{n\in F}\frac{1}{n^2t_n}t_n=
\sum_{n\in F}\frac{1}{n^2}\leq\sum_{n\in \mathbb{N}}\frac{1}{n^2}<\infty,$$
i.e. $h\in L_1(M,\tau)$.

Put $\mu(\cdot)=\tau(h\cdot)$. Since $y\in M_f$, it follows that $y\in L_1(M,\mu)$. Therefore $\mu(y)<\infty$.

On the other hand
$$hy=\sum_{n\in F}\frac{1}{n^2t_n}p_n \sum_{n\in \mathbb{N}}n^2p_n=\sum_{n\in F}\frac{1}{t_n}p_n,$$
and thus
$$\mu(y)=\tau(hy)=\sum_{n\in F}\frac{1}{t_n}\tau(p_n)=\sum_{n\in F}\frac{1}{t_n}t_n=\sum_{n\in F}1=|F|,$$
where $|F|$ is the cardinality of the set $F$.
Since $\mu(y)<\infty$ this implies that $F$ is a finite set. Let $k=\max\{n:n\in F\}$. Then $\tau(p_n)=0$ for all $n>k$,
and since $\tau$ is faithful we have that $p_n=0$ for all $n>k$, i.e.
$e_{(n+1)^2}=e_{n^2}$. But $e_{n^2}\uparrow\textbf{1}$ and thus
$e_{n^2}=\textbf{1}$ for all $n>k$.
This means that $0\leq x\leq (k+1)^2\textbf{1}$, i.e. $x\in Z$.

The proof is complete.  $\blacksquare$

\begin{pr}
Let $M$ be a type $I_n,\, n\in\mathbb{N}$ von Neumann algebra. Then $M_f=M.$
\end{pr}

\emph{Proof.} By \cite[Ch. V, Theorem 1.27]{Tak} the von Neumann algebra $M$
of type $I_n\,\,\,(n\in\mathbb{N})$ can
be represented as $M=Z\otimes B(H_n),$ where $Z$ is the center  $M$ and
$H_n$ is the $n$-dimensional Hilbert space.
Put $\mathcal{F}_Z=\{\tau|_{Z}:\tau\in\mathcal{F}\}.$
 Therefore from Proposition 3.1 we obtain

$$M_f=\bigcap\limits_{p\in[1,\infty)}\bigcap\limits_{\tau\in \mathcal{F}}L_{p}(M,\tau)=
\bigcap\limits_{p\in[1,\infty)}\bigcap\limits_{\mu\in \mathcal{F}_Z}
L_{p}(Z,\mu)\otimes B(H_n)=$$
$$=\left(\bigcap\limits_{p\in[1,\infty)}\bigcap\limits_{\mu\in \mathcal{F}_Z}L_{p}(Z,\mu)\right)\otimes B(H_n)=
Z_{f}\otimes B(H_n)=$$
$$=Z\otimes B(H_n)=M,$$
i.e. $M_f=M.$

The proof is complete.  $\blacksquare$

\begin{pr}
Let $M$ be a finite von Neumann algebra which is isomorphic to the direct sum of an infinite number of
homogeneous type $I_n\,\,\,(n\in\mathbb{N})$ von Neunamm algebras. Then $M_f\neq M.$
\end{pr}

\emph{Proof.} Suppose that $M={\sum\limits_{k\in K}}^{\oplus}M_k,$ where  $K$ is an infinite subset
of $\mathbb{N},$ and $M_k$ is a homogeneous type $I_k$ von Neumann algebra.

Since the set $K$ is infinite, there exists a sequence $\{k_n\}\subset K$
such that $k_n\geq 2^n$ for all $n\in \mathbb{N}$.
We have that
$$M_{k_n}=Z_{k_n}\otimes B(H_{k_n}),$$
where $Z_{k_n}$ is the center of $M_{k_n}$ and
$$N_n=\textbf{1}_n\otimes B(H_{2^n})\subset M_{k_n}.$$
Therefore the algebra $M$ contains a subalgebra *-isomorphic
 to the algebra $N={\sum\limits_{n\in \mathbb{N}}}^\oplus N_n$.

Hence, without loss of  generality we may assume that  $M={\sum\limits_{n\in \mathbb{N}}}^{\oplus}N_n,$
where $N_n=B(H_{2^{n}})$~  is  the algebra of all $2^n\times 2^n$ matrices over $\mathbb{C}.$ On each $N_n$ we consider
the unique tracial state (i.~e. normalized f.n.f. trace) $\mu_n$ and define on $M$ the following f.n.f. trace
$$\tau(x)=\sum\limits_{n\in \mathbb{N}}2^{-n}\mu_{n}(x_{n}),$$
where $x={\sum\limits_{n\in \mathbb{N}}}^{\oplus}x_n\in M.$ Then every f.n.f. trace $\mu$ on $M$ has the form
$$\mu(x)=\tau(hx)=\sum\limits_{n\in \mathbb{N}}2^{-n}\mu_{n}(h_{n}x_{n})=
\sum\limits_{n\in \mathbb{N}}2^{-n}\alpha_{n}\mu_{n}(x_{n}),$$
where
$$h={\sum\limits_{n\in \mathbb{N}}}^{\oplus}h_n=
{\sum\limits_{n\in \mathbb{N}}}^{\oplus}\alpha_n\textbf{1}_{n}\in L_1(M,\tau),
$$
i.~e. $\alpha_n>0,\,\,n\in \mathbb{N},$$\sum\limits_{n\in \mathbb{N}}2^{-n}\alpha_n<\infty.$

Take a minimal projection $p_n$ in each $N_n=B(H_{2^{n}}).$ Then $\mu_n(p_n)=\frac{\textstyle 1}{\textstyle 2^n}.$

Consider the unbounded element  $x={\sum\limits_{n\in \mathbb{N}}}^{\oplus}n p_n$ in $S(M)\setminus M$ and
let us prove that $x\in M_f.$  For every f.n.f. trace $\mu$ on $M$ one has that
$$\mu(x^p)=\sum\limits_{n\in \mathbb{N}}2^{-n}\alpha_n\mu_n(n^p p_n)=
\sum\limits_{n\in \mathbb{N}}2^{-n}\alpha_nn^p2^{-n}<\infty,$$
because $n^p2^{-n}<1$ for sufficiently large $n\in\mathbb{N}.$ Therefore $x\in L_p(M,\mu)$
for all $p\geq 1$ and every f.n.f. trace $\mu\in\mathcal{F},$ i.~e. $x\in M_f.$

The proof is complete.  $\blacksquare$

\begin{pr}
Let $M$ be  a type $II_1$ von Neumann algebra with a f.n.f. trace $\tau.$ Then $M_f\neq M$.
\end{pr}

\emph{Proof.} Suppose that the trace $\tau$ is normalized, i.~e.
$\tau(\textbf{1})=1$, and denote  by $\Phi$ the canonical
center-valued trace on $M$. Since $M$ is of type $II_1$ there exists a
projection $p_1$ such that
$$p_1\sim \textbf{1}-p_1.$$
Therefore from $\Phi(p_1)+\Phi(p_1^{\perp})=\Phi(\textbf{1})=\textbf{1}$ and
$\Phi(p_1)=\Phi(p_2)$ we obtain that
$$\Phi(p_1)=\Phi(p_1^{\perp})=\frac{\textstyle 1}{\textstyle 2}\textbf{1}.$$

Suppose that we have constructed mutually orthogonal projections
$p_1,\,p_2,\,
\cdots, \,p_n$ in $M$ such that
$$\Phi(p_k)=\frac{\textstyle 1}{\textstyle 2^{k}}\textbf{1},\,k=\overline{1, n}.$$
Set $e_n=\sum\limits_{k=1}^{n}p_k.$ Then
$\Phi(e_n^{\perp})=\frac{\textstyle 1}{\textstyle 2^{n}}\textbf{1}.$
Now take a projection $p_{n+1}\leq e_n^{\perp}$ such that
$$p_{n+1}\sim e_n^{\perp}-p_{n+1},$$
i.~e.
$$\Phi(p_{n+1})=\frac{\textstyle 1}{\textstyle
2^{n+1}}\textbf{1}.$$

In this manner we obtain a sequence $\{p_n\}_{n\in\mathbb{N}}$  of mutually orthogonal
projections such that
$$\Phi(p_n)=\frac{\textstyle
1}{\textstyle 2^{n}}\textbf{1},\,n\in\mathbb{N}.$$
It is clear that  $\tau(p_n)=\tau(\Phi(p_n))=\frac{\textstyle 1}{\textstyle 2^{n}}, \,\,\,n\in\mathbb{N}.$

From
$$\sum\limits_{n=1}^{\infty}||np_n||_{1}^{\tau}=
\sum\limits_{n=1}^{\infty}\tau(np_n)=
\sum\limits_{n=1}^{\infty}\frac{\textstyle n}{\textstyle
2^{n}}<\infty,$$
it follows  that the element  $x=\sum\limits_{n=1}^{\infty}np_n$ belongs
to $ L_{1}(M, \tau),$ and it is unbounded, i.~e.  $ x\notin M.$

On the other hand for an arbitrary central element $h\in L_{1}(M, \tau), h> 0,$  and
 $n\in\mathbb{N}$  we have
$$\tau(hp_n)=\tau(\Phi(hp_n))=\tau(h\Phi(p_n))=\tau(h\frac{\textstyle
1}{\textstyle 2^{n}}\textbf{1})=\frac{\textstyle 1}{\textstyle
2^{n}}\tau(h).$$
Therefore for an arbitrary f.n.f. trace $\mu$ on $M$ with  $\frac{\textstyle d\mu}{\textstyle d\tau}=h$ we have
$$\mu(|x|^{p})=\mu(x^{p})=\tau(h x^{p})=
\tau(h\sum\limits_{n=1}^{\infty}n^{p}p_n)=$$
$$=\sum\limits_{n=1}^{\infty}n^{p}\tau(hp_n)=\tau(h)\sum\limits_{n=1}^{\infty}\frac{\textstyle
n^{p}}{\textstyle 2^{n}}<\infty,$$
 i.~e. $x\in L_p(M,\mu)$ for all $p\geq 1$ and every f.n.f. trace $\mu.$
 Therefore $x\in M_f\setminus M.$

 The proof is complete.  $\blacksquare$

 \emph{Proof of Theorem 3.2.} The implication $(i)\Rightarrow (ii)$
follows from Propositions 3.3 and 3.4, while $(ii)\Rightarrow (i)$
follows from Propositions 3.2.

The proof is complete.  $\blacksquare$

Now let us describe continuous linear functionals on the space $(M_f, t)$.

\begin{thm}
Given any $\mu\in \mathcal{F}$, $1<q<\infty$, and $a\in L_q(M,\mu)$ the functional $\varphi(x)=\mu(xa)$, $x\in M_f$, is a
continuous linear functional on $(M_f, t)$. Conversely for any continuous linear functional
 $\varphi$ on $(M_f,t)$ there exist
$\mu\in\mathcal{F},\,\, 1<q<\infty,\,\, a\in L_{q}(M,\mu)$ such that
$$\varphi(x)=\mu(x a),\,\,\,x\in M_f.$$
\end{thm}

\emph{Proof.} Let   $\mu\in\mathcal{F},\,\, 1<q<\infty,\,\, a\in L_{q}(M,\mu).$ Put
$$\varphi_a(x)=\mu (x a),\,\,\,x\in M_f.$$
Take  $p\in\mathbb{R}$ such that
$\frac{\textstyle 1}{\textstyle p}+\frac{\textstyle 1}{\textstyle q}=1.$
Since
$$|\varphi_a(x)|=|\mu (x a)|\leq ||a||_{q}^{\mu}||x||_{p}^{\mu}$$
for all $x\in M_f,$ one has that $\varphi_a$ is a continuous linear functional on   $(M_f,t).$

Conversely, let $\varphi$ be a continuous linear functional on $(M_f,t).$
 By \cite[Corollary 1 on p.43]{Yo} there exist $\mu\in\mathcal{F},\,\, 1\leq p<\infty,\,\,c>0$,
 such that
 $$|\varphi(x)|\leq c||x||_{p}^{\mu}$$
for all $x\in M_f.$
 Since $M\subset M_f$ and $M$ is $\|\cdot\|_p^\mu$-dense in $L_p(M,\mu)$,
  the functional  $\varphi$ can be uniquely extended onto $L_p(M, \mu)$.
  By \cite[Theorem 4.4]{Yed} there exists $a\in L_{q}(M,\mu),\,
\frac{\textstyle 1}{\textstyle p}+\frac{\textstyle 1}{\textstyle q}=1$,
such that $$\varphi(x)=\mu(x a)$$ for all
$x\in L_{p}(M,\mu).$ In particular
$$\varphi(x)=\mu (x a)$$ for all
$x\in M_f,$ i.e. $\varphi=\varphi_a.$

 The proof is complete.  $\blacksquare$

If the von Neumann algebra $M$ is a factor
 then it has  a unique (up to a scalar multiple) f.n.f. trace $\mu.$
In this case the finite tracial algebra $M_f$ coincides with the Arens
 algebra $L^{\omega}(M,\mu)$ and the topology
$t$ merges to the topology $t^{\mu}$ generated by  the system
of norms $\{\|\cdot\|^{\mu}_{p}\}_{p\geq 1}.$
The following theorem describes the general case where this phenomenon occurs.

Recall some notions from the theory of linear topological spaces.
Let $E$ be a locally convex linear topological space.
An absolutely convex absorbing set in $E$ is called a barrel.
If each barrel in $E$ is a neighborhood of zero, then $E$ is
said to be $a$ \emph{barreled space}.

It is known  (\cite{Yo}, Theorem 2, p.200 ) that every reflexive locally convex space is barreled.

\begin{thm}
Let $M$ be a finite von Neumann algebra and suppose that $\mathcal{F}\neq \varnothing$ is the family
of all f.n.f. traces on $M.$ The following conditions are equivalent:

(i)   $M_f=L^{\omega}(M,\mu)$ for some (and hence for all)  $\mu\in \mathcal{F};$

(ii)  $(M_f,t)$ is metrizable;

(iii) $(M_f; t)$ is reflexive;

(iv) the center $Z$ of $M$ is finite dimensional, i.~e.
$M=\sum\limits_{i=1}^{m}M_{i},$ where all $M_i$  are  $I_n$-factors or $II_1$-factors.
\end{thm}

\emph{Proof.} Suppose that $Z$ is finite dimensional. Then $M$ is a finite direct
sum of factors $M_i,\,\,\,i=\overline{1,k}.$ Then  for each factor $M_i$
the algebras $(M_i)_f$ and $L^\omega(M_i, \mu_i)$ coincide and the topology $t_i$ is the same as $t_{i}^{\mu_i}.$
Therefore
$$M_{f}=(\sum\limits_{i=1}^{n}M_{i})_{f}=\sum\limits_{i=1}^{n}(M_{i})_{f}=
\sum\limits_{i=1}^{n}L^\omega(M_i, \mu_i)=L^\omega(M, \mu),$$
where $\mu={\sum\limits_{i=1}^{n}}\mu_{i}\in \mathcal{F},$ i.~e.
 $M_{f}=L^\omega(M, \mu).$

Now since the topology $t^{\mu}$ on the Arens algebra
$L^\omega(M, \mu)$ is metrizable \cite{Abd} it follows that $t=t^\mu$  is also metrizable.

It is known \cite{Abd2} that for finite traces $\mu$ the Arens algebra $(L^\omega (M, \mu), t^\mu)$ is reflexive  and hence $(M_f, t)$ is also reflexive.

Therefore $($\emph{iv}$)$ implies $($\emph{i}$)$, $($\emph{ii}$)$ and $($\emph{iii}$)$.

$($\emph{i}$)$ $\Rightarrow$ $($\emph{iv}$)$. Suppose
that $M_f=L^\omega(M, \mu)$ for an appropriate
 $\mu\in \mathcal{F}$. Then there exists a sequence
of mutually orthogononal projections $\{p_n\}$ in $Z$
such that $p_n\neq 0$ for all $n\in \mathbb{N}$. Since the trace $\mu$ is finite one has
that $\sum\limits_{k=1}^{\infty}\mu(p_k)<\infty$ and hence there
is a subsequence $\{n_k: k\in \mathbb{N}\}$ such
that $\mu(p_{n_k})\leq \frac{\textstyle 1}{\textstyle  2^k}$ for all $k$.

Set
$$x=\sum\limits_{k=1}^\infty k p_k$$
For $p\geq 1$ we have
$$\mu(|x|^p)=\sum\limits_{k=1}^\infty k^p\mu(p_k)\leq\sum\limits_{k=1}^\infty k^p\frac{1}{2^k}<\infty,$$
and hence $x\in L^\omega(M, \mu)=M_f$.

On the other hand $x$ is a central element in $M_f$ and Proposition 3.1
implies that $x\in Z(M_f)=Z\subset M$. But it is clear that the
element $x$ is unbounded, i.e. $x\notin M$. The contradiction shows that $Z$ is finite dimensional.

$(ii)\Rightarrow (iv).$ Suppose that $(M_f,t)$ is metrizable.
By Theorem 3.1 it is complete and hence it is a Fre$^{\prime}$chet
space. In particular the center of $M_f$ which coincides with $Z_f$ is also a Fre$^{\prime}$chet space.
By  Proposition 3.1
$Z_f=Z$ and hence $Z$ is a Fre$^{\prime}$chet space with respect to the induced topology $t_z=t|_{Z}.$

Consider the identity mapping
$$I:(Z, \|\cdot\|_\infty)\rightarrow (Z,t_z)$$
where $\|\cdot\|_\infty$ is the operator norm on $Z.$
From the inequalities
$$
\|x\|_{p}^{\mu}\leq C_{p}^{\mu}\|x\|_\infty
$$
(where $C_{p}^{\mu}$ is  an appropriate constant for each $p\geq1,$
 $\mu\in\mathcal{F}$) it follows that the mapping $I$ is continuous.
Since $(Z,t_z)$ is a Fre$^{\prime}$chet space, from Banach theorem on the inverse operator
(\cite{Yo}, Chapter II, Section 5) we obtain that the inverse mapping
$$I^{-1}:(Z,t_z)\rightarrow (Z, \|\cdot\|_\infty)$$
is also continuous. This means that for some $p\in[1,\infty)$ and an appropriate
$\mu\in \mathcal{F}$ there exists  a constant $K_{p}^{\mu}$ such that
\begin{equation}
\|x\|_\infty\leq K_{p}^{\mu}\|x\|_{p}^{\mu}
\end{equation}
for all $x\in Z$ (\cite{Yo}, Theorem 1, p. 42).

Now suppose that $\dim Z=\infty$. There exists a sequence $\{p_n\}$ of projections in $Z$
 such that $p_n\uparrow \textbf{1}$, $p_n\neq p_{n+1}$.
Thus $p_n^\perp\neq 0$, $\tau(p_n^\perp)\rightarrow 0$, i.e.
 $\|p_n^\perp\|_p^\mu\rightarrow 0$. From the inequality (1) we obtain that
$\|p_n^\perp\|_\infty \rightarrow 0$.

On the other hand $\|p_n^\perp\|_\infty =1$. This contradiction implies that $Z$ is finite dimensional.

$(iii)\Rightarrow(iv)$. Suppose that $M_f$ is reflexive.
 Then the center $Z(M_f)=Z$ is also reflexive as a closed subspace of
a reflexive space.

The set
$$B=\{x\in Z: ||x||_\infty\leq 1\}$$
is a barrel in $(Z, t)$ and since $Z$ is reflexive, we have that $B$ is a neighborhood
 of zero in $Z$. Therefore there exist $p\geq 1$,
$\mu \in \mathcal{F}$ and $\varepsilon >0$ such that
$$\{x\in Z: \|x\|_p^\mu\leq \varepsilon\}\subseteq B$$
i.e.
$$\|x\|_\infty\leq \varepsilon^{-1}\|x\|_p^\mu$$
for all $x\in Z$. From this as above it follows that $Z$ is finite dimensional.

The proof is complete.  $\blacksquare$

\textbf{Remark.} In the von Neumann algebra $M$ the operator topology is stronger than the topology  $t$,
$t$ is stronger than    $t^{\mu}$, and   $t^{\mu}$ is stronger than each $L_p$-norm topology for any $p\geq 1$.

\section{Derivations on Finite Tracial Algebras}

Derivations on unbounded operator algebras, in particular on
various algebras of measurable operators affiliated with von
Neumann algebras, appear to be a very attractive special case of
 general  unbounded derivations on operator algebras.

Let  $A$ be an algebra over the complex number. A linear operator
$D:A\rightarrow A$ is called a \emph{derivation} if it satisfies the
identity  $D(xy)=D(x)y+xD(y)$ for all  $x, y\in A$ (Leibniz rule).
Each element  $a\in A$ defines a derivation  $D_a$ on $A$ given as
$D_a(x)=ax-xa,\,x\in A.$ Such derivations $D_a$ are said to be
\emph{inner derivations}.

In \cite{Alb2} we have investigated and completely
described derivations on the algebra $LS(M)$ of all
locally measurable operators affiliated with a type I von Neumann
algebra $M$ and on its various subalgebras. Recently the
above conjecture was also confirmed for the type I case in the
paper \cite{Ber1} by a representation of measurable operators as
operator valued functions. Another approach to similar problems in
$AW ^{*}$-algebras of type I was suggested in the recent paper
\cite{Gut}.

In the paper \cite{Alb} we have proved the spatiality of
derivations on the non commutative Arens algebra
$L^{\omega}(M, \tau)$ associated with an arbitrary von
Neumann algebra $M$ and a faithful normal semi-finite
trace $\tau.$ Moreover if the trace $\tau$ is finite then every
derivation on $L^{\omega}(M, \tau)$ is inner.

In this section we prove that each derivation on a finite tracial algebra is inner.

The following result is an immediate corollary of \cite [Proposition 3.6]{Ayu}.

\begin{lm} Let $M$ be a von Neumann algebra with a faithful
 normal trace $\tau$. Given any derivation $D:M\rightarrow L^\omega (M, \tau)$
there exists an element $a\in L^\omega(M, \tau)$ such that
$$D(x)=ax-xa, \,\,\, x\in M.$$
\end{lm}

Further we need also the following assertion from \cite[Proposition 6.17]{Ber1}.

\begin{lm} Let $A$ be a *-subalgebra of $LS(M)$ such that $M\subseteq A$
and $A$ is solid (that is, if $x\in LS(M)$ and $y\in A$
satisfy $|x|\leq |y|$ then $x\in A$). If $\omega \in LS(M)$ is such that
 $[\omega, x]\in A$ for all $x\in A$, then there exists $\omega_1\in A$
such that $[\omega, x]=[\omega_1, x]$ for all $x\in A$.
\end{lm}

The main result of this section is the following theorem.

\begin{thm}  Let $M$ be a von Neumann algebra with
 a faithful normal finite trace $\tau$.
  If $A\subseteq L^\omega(M, \tau)$ is a solid
*-subalgebra such that $M\subseteq A$,
 then every derivation on $A$ is inver.
\end{thm}

\emph{Proof.} Since
$A\subseteq L^\omega(M, \tau)$, by Lemma 4.1 there
 exits an element $a\in L^\omega(M, \tau)$ such that
\begin{equation}
D(x)=ax-xa, \,\,\,\,\,\,\,\, x\in M.
\end{equation}

Let us show that in fact
\begin{equation}
D(x)=ax-xa,\,\,\,\mbox{for all} \,\,\,\, x\in A.
\end{equation}
Consider $x\in A, \,\, x\geq 0$. Then $(1+x)^{-1}\in M$. From the
 Leibniz rule it follows that for each invertible $b\in A$ one has
$$D(b)=-b D(b^{-1}) b.$$
Therefore
$$D(x)=D(\mathbf{1}+x)=-(\mathbf{1}+x)D((\mathbf{1}+x)^{-1})(\mathbf{1}+x).$$
On the other hand since $(1+x)^{-1}\in M$ the equality (2) implies that
$$D((1+x)^{-1})=a(\mathbf{1}+x)^{-1}-(\mathbf{1}+x)^{-1}a.$$
Therefore
$$-(\mathbf{1}+x)D((\mathbf{1}+x)^{-1})(\mathbf{1}+x)=
-(\mathbf{1}+x)[a(\mathbf{1}+x)^{-1}-(\mathbf{1}+x)^{-1}a](\mathbf{1}+x)=$$
$$=-(\mathbf{1}+x)a+a(\mathbf{1}+x)=ax-xa,$$
i.e.
$$D(x)=ax-xa,\,\,\,x\in A, \,\,\, x\geq 0.$$

Since each element from $A$ is a finite linear combination
 of positive elements, we obtain the equality (3) for arbitrary $x\in A$.

Now since $A$ is a solid *-subalgebra in $L^\omega(M, \tau)$ containing $A$,
Lemma 4.2 implies that the element $a$ implementing
the derivation $D$ may be chosed from the algebra $A$, i.e.
$$D(x)=ax-xa,\,\,\,x\in A$$
for an appropriate $a\in A.$

The proof is complete.   $\blacksquare$

Since the algebra $M_f$ is a solid *-subalgebra of
$L^\omega(M, \tau)$ and contains $M$, we obtain the following result.

\begin{cor}
If $M$ is a von Neumann algebra with a faithful normal trace, then every derivation on $M_f$ is inner.
\end{cor}

\vspace{1cm}

\textbf{Acknowledgments.} \emph{Part of this work was done within
the framework of the Associateship Scheme of the Abdus Salam
International Centre  for Theoretical Physics (ICTP), Trieste,
Italy. The first author would like to thank ICTP for the kind hospitality and for providing financial support
 and all facilities (July-August, 2009).This work is supported in part by the DFG 436 USB 113/10/0-1
project (Germany).}

\end{document}